\documentclass{article}
\usepackage[utf8]{inputenc}
\usepackage{xcolor}
\usepackage{lscape}
\usepackage{array}
\usepackage{hyperref}

\usepackage[margin=1.5in]{geometry}
\usepackage{comment}
\usepackage[inline]{enumitem}
\usepackage{amsthm}

\usepackage{soul}

\usepackage[colorinlistoftodos]{todonotes}

\newcommand{\ron}[1]{\todo[size=\small,inline,color=magenta!50]{#1 \\ \hfill --- Ron}}

\newcommand{\rachel}[1]{\todo[size=\small,inline,color=violet!30]{#1 \\ \hfill --- Rachel}}

\definecolor{munsell}{rgb}{0.0, 0.5, 0.69}

\title{On definitions of ``mathematician"}

\author{Ron Buckmire
\and Carrie Diaz Eaton \and Joseph E. Hibdon, Jr. \and Katherine M. Kinnaird \and Drew Lewis \and Jessica Libertini \and Omayra Ortega \and Rachel Roca \and Andr\'es R. Vindas-Mel\'endez}

\date{\today}
\begin{document}

\maketitle




\begin{abstract}

The definition of who is or what makes a ``mathematician" is an important and urgent issue to be addressed in the mathematics community. 
Too often, a narrower definition of who is considered a mathematician (and what is considered mathematics) is used to exclude people from the discipline---both explicitly and implicitly. 
However, using a narrow definition of a mathematician allows us to examine and challenge systemic barriers that exist in certain spaces of the community. 
This paper explores and illuminates tensions between narrow and broad definitions and how they can be used to promote both inclusion and exclusion simultaneously. 
In this article, we present a framework of definitions based on identity, function, and qualification and exploring several different meanings of {\em mathematician}.
By interrogating various definitions, we highlight their risks and opportunities, with an emphasis on implications for broadening and/or narrowing participation of underrepresented groups.
\end{abstract}

As we reflect on our mathematical communities' journeys towards social justice, we have been thinking about how to use mathematics and data to keep these communities accountable to and focused on this goal; specifically, we are interested in data about who is included, who is excluded, who is overrepresented, and who is underrepresented. 
As we began shaping possible research projects, we entered into conversations about data availability and data collection, and a natural question emerged: We began to ask, ``How do we define the mathematical community?" 

Thinking about definitions is a common first act of mathematical inquiry. 
Thinking about how we are defining ``mathematical community" led us to ask for a specific definition for ``mathematician." 
The co-authors reflected on how designating who is and who is not a mathematician is itself a complex social justice issue. 
It is difficult to consider the question of defining a mathematician without also considering the question of what is mathematics.  
Of course, one can reasonably define mathematics in numerous ways. 

For the purposes of this paper, we will use ``mathematics'' as an umbrella term for the mathematical sciences, broadly construed.
Additionally, while we, the authors, believe the ``mathematics community" is the set of all people who identify as mathematicians, throughout the paper, we will often use the term ``mathematics community'' or ``mathematicians'' to mean the set of mathematicians that results from explicitly employing a particular definition. 
We understand that a ``set of mathematicians'' does not imply ``community,'' but use this word to intentionally humanize our language.
This approach allows us to explore the impact of different definitions of mathematician on the make-up of the community that results when a given definition is used to form the boundary between the in-group and out-group.

This paper seeks to explicitly uncover and unpack assumptions and ideologies associated with the term ``mathematician." 
This paper will do three things: 1) provide a framework based explicitly on three different ways to define the term ``mathematician" (identity, function and qualification) informed by the co-authors' lived experiences that complements existing, broader mathematical identity frameworks found in the literature; 2) analyze and discuss the definitions, assumptions and implications of the term ``mathematician";  and 3) act as an entry point for readers of this paper to engage in self-reflection about their identities as mathematicians and the importance and ramifications of different definitions of the term. The goal is to provide a framework for others who work with data on the ``mathematics" community to clarify their definitions and the reasoning for their definition choice for the specific issue they are addressing.

We note that we write this paper from the perspective of people grounded in various mathematical backgrounds and with different lived experiences, rather than from a grounding in social science research around identity. In Section~\ref{subsec:positionality}, we provide more detailed positionality statements.

The rest of the paper is organized as follows. We conclude this section with a brief discussion of the relevant literature related to mathematics identity. In Section~\ref{sec:defs} we provide our definitional framework for ``mathematician" based on three aspects of the term: qualification, identity, and function, along with broad and narrow examples of each definition-type. This is followed by our analysis of the implications of our definitional framework in  Section~\ref{sec:analysis} that begins with our positionality statements and includes  reflections on them and features a visualization of how the co-authors reflect our definitional framework given in Table~\ref{table:position}. This section also includes our implementation of a polarity analysis that discusses the risks and opportunities associated with the broad and narrow versions of each of the definition-types in our ``mathematician" definitional framework along with a visualization of the analysis in Figure~\ref{fig:summary}. The paper ends with two short sections with conclusions in Section~\ref{sec:conclusions} and discussion of future work in section\ref{sec:future}.

\subsection*{Review of Mathematics Identity Scholarship}
The question of defining the term ``mathematician'' is intimately related to the idea of {\em mathematics identity}. 
Mathematics identity is well-studied in the mathematics education literature, particularly in K-12 education and with respect to the formation of mathematical identity, while the definition of ``mathematician'' is often treated as an axiom or universally accepted given.  
Our paper builds on the mathematics identity scholarship, analyzing these identity constructs to define ``mathematician" with a focus on practicing mathematicians. 

In an extensive literature review of the topic, Darragh \cite{darragh} broadly categorizes definitions of mathematics identity into five groups: participative, narrative, discursive, psychoanalytic, and performative.  
A participative definition is one that defines identity through participation and engagement with a social group. 
A narrative approach defines identity through the stories we tell about mathematics. 
Sfard and Prusak \cite{SfardPrusak} go as far as to say that identities are stories. 
A discursive perspective of identity views it through the language people use when they interact.
This can be viewed either as the language people use when they do mathematics, or to refer to broader societal narratives about mathematics.
A psychoanalytic definition describes identity in terms of interactions between conscious and unconscious processes.  
A performative conceptualization holds that identity is constituted through one's actions, repeated over time. 

Black et al. \cite{black} use three categories of definitions. 
The first two are discursive and psychoanalytic, as above, while they term their third category as sociocultural, which is a broadening of the participative definition above.

Gee \cite{Gee} outlines four overlapping ways to view identity.  
A nature identity captures an intrinsic state of being.
An institutional identity is where identity  granted by an institution.
Discourse identities are, how one discourses and engages with others, and how others discourse about one.  
Affinity identities are centered around communities of shared identity to which one chooses to belong, similar to Darragh's participative definition.

Darragh, following Guti\'errez \cite{gutierrez} and Gee \cite{Gee}, further highlights a dichotomy between mathematics identity as something intrinsic that an individual has acquired, and mathematics identity as a set of actions, i.e., things that an individual does.

Guti\'errez \cite{gutierrez} highlights how (mathematics) identity and power are intertwined.   Being granted access to the community of ``mathematicians," however defined, grants one access to these power structures.  
We now turn our attention to definitions of ``mathematician,'' building upon the above literature on mathematics identity and power.
We then analyze the opportunities and risks of various definition choices will integrate these perspectives to focus both on on implications for identity formation, maintenance, and power.

\section{Definitions}\label{sec:defs}

In this section, we present our definitional framework for the term ``mathematician." 
This framework is based on three kinds of definitions, each of which can be tuned across a spectrum. 
These definitions with their accompanying ranges provide opportunities for a wide range of insights in exploring one's identity from multiple perspectives. 
Building from these definitions, we produce another set of  definitions by constructing unions and intersections of the first three types of definitions.

From our lived experiences and discussions, we recognized the need to explicitly design a framework that more specifically serves to define the ``mathematical community" and a ``mathematician." 
This framework below---highlighting function, identity, and qualification---is not meant to contradict previous work done on identity, but to complement and focus it on the unique structures and politics surrounding mathematics.
Our framework, like others, explicitly recognizes complexities in the social nature of identity, while addressing specific issues the mathematics community faces.

We further present a variety of definition types for who is and what makes a mathematician. 
Of course, this list is not exhaustive. 
For each type---function, identity, and qualification---we offer three example definitions, including narrow and broad example definitions for mathematicians as well as a ``hybrid'' example that is somewhere between the broad and narrow extremes. 
Recall that we defined ``mathematics'' in the context of this paper as the mathematical sciences, but we note that the way mathematics is defined could also broaden or narrow the pool of those that satisfy the definitions given below.
We emphasize that these are examples of function-based, qualifications-based, and identity-based definitions that we, as the authors, have developed. 
We encourage others to explore other definitions that also align with the given framework. 

\subsection{Function-based Definitions}

A function-based definition for a mathematician is a definition that hinges on what they do (or how they ``function'') on a daily basis.
Function-based examples may include a collection of specific actions or a more general list of habits or holistic activities. 

Below we provide both broad and narrow examples of function-based definitions followed by a hybrid example definition, which we emphasize are our examples of function-based definitions that we, as authors, have developed.

\newtheorem{functional_ex_definition}{\emph{Function-based  Definition}}

\begin{functional_ex_definition}[Broad]\label{def:functional-broad}
A mathematician is a person who uses mathematical concepts, tools or techniques to study and solve problems.
\end{functional_ex_definition}

Function-based Definition \ref{def:functional-broad} may include all humans. 
If counting is considered a mathematical technique, then anyone who counts for a problem-solving purpose can then be defined as a mathematician. 
This definition is often invoked in books such as Keith Devlin's \textit{The Math Instinct} \cite{devlin2009} in order to make the case that ``we are all mathematicians.''
This use is related to a broader effort to boost confidence in more people to empower them to be successful in their mathematical work and to combat the prevalent and pernicious idea that ``Some people are not 'math' people."

\begin{functional_ex_definition}[Narrow]\label{def:functional-narrow}
A mathematician is a person who proves theorems using proof techniques.
\end{functional_ex_definition}

Function-based Definition \ref{def:functional-narrow} is often implied in the construction of mathematics majors.
It manifests in, and is reinforced by, the requirement of a proofing course as the foundational ``methods'' courses in the curricular design of the undergraduate mathematics major at some institutions.

\begin{functional_ex_definition}[Hybrid]\label{def:functional-hybrid}
A mathematician is a person who, as part of their daily work, employs mathematical techniques and tools to solve mathematical and other problems.
\end{functional_ex_definition}

There are many function-based example definitions that fall between the narrow and broad example definitions. 
However, Function-based Definition \ref{def:functional-hybrid} might be used to identify jobs held by such mathematicians, in order to advise students about career options beyond academia.

\subsection{Qualifications-Based Definitions}
In a qualification-based definition, a mathematician is someone who has cleared a normative bar established to certify knowledge or accomplishment.
This bar or ``qualification" is something that can be moved depending on what someone is trying to measure. 
The construction of this definition implies the existence of an external body that certifies that the qualification has occurred. 

\newtheorem{qualification_ex_definition}{\emph{Qualification-Based Definition}}
\begin{qualification_ex_definition}[Broad]\label{def:qualification-broad}

    A mathematician is someone who has completed a course in mathematics at any level.
\end{qualification_ex_definition}

In Qualification-Based Definition~\ref{def:qualification-broad} the external body that has certified the qualification (the ``completion'') would be the teacher of the course. 
We could be more specific in this definition, by choosing an example like ``being able to count to ten" or ``can correctly multiply two-digit integers" or ``has taken a course in a specific mathematics topic like trigonometry."

\begin{qualification_ex_definition}[Narrow]\label{def:qualification-narrow}
    A mathematician is someone who has a Ph.D. in mathematics and currently holds a research position in mathematics. 
\end{qualification_ex_definition}

Qualification-based Definition \ref{def:qualification-narrow} is often operationalized in order to make hiring decisions consistent with terminal degree expectations for post-secondary accreditation review.

\begin{qualification_ex_definition}[Hybrid]\label{def:qualification-hybrid}
    A mathematician is someone who holds a Bachelor's degree or higher in mathematics or mathematics coursework, where the mathematics coursework allows for one to teach secondary education mathematics. 
\end{qualification_ex_definition}

The teacher-credentialing process in mathematics uses a qualifications-based definition which equates a set of post-secondary coursework with sufficient pedagogical content knowledge for teaching. Guidelines of these criteria are given by the National Council of Teachers of Mathematics \cite{NCTM}. 
However the Qualification-Hybrid definition is a moving definition based on some arbitrary criteria. 
Many of these qualifications have been set initially, but not  updated enough to stay in line with current societal needs.

\subsection{Identity-based Definitions}
Unlike the qualifications-based and the function-based definitions which lie on a narrow-broad spectrum, we can consider identity-based definitions along an internal-external spectrum.
An internal identity-based definition centers the individual perspective, while the external identity-based definition relies on the collective perception of others on the individual's mathematician status. 

These identity-based example definitions highlight the overlap between Gee's affinity-identity and discursive-identity \cite{Gee}.

\newtheorem{identity_ex_definition}{\emph{Identity-Based Definition}}
\begin{identity_ex_definition}[Internal]\label{def:identity-internal}
    A mathematician is anyone who says they are a mathematician.
\end{identity_ex_definition}
Bottom line, if you decide you are a mathematician, by Identity-Based Definition \ref{def:identity-internal}, you are. If you decide you are not, then you are not. This definition invokes the idea of self-determination of identity, which can be a particularly powerful activity to those who experience oppression because of their identities \cite{morton}. 

\begin{identity_ex_definition}[External]\label{def:identity-external}
    A mathematician is anyone whom others say is a mathematician.
\end{identity_ex_definition}

Identity-Based Definition \ref{def:identity-external} requires others to invoke their identity as a mathematician, and then requires enough consensus so that the title is formally or informally bestowed. Formal versions of this type of identity validation could include earning an award given only to mathematicians.

\begin{identity_ex_definition}[Hybrid]\label{def:identity-broad}
    A mathematician is anyone who someone else says is a mathematician.
\end{identity_ex_definition}
Like Identity-Based Definition \ref{def:identity-external}, Identity-Based Definition \ref{def:identity-broad} is still external, but only requires one person's criteria for mathematician is met. 
Therefore, it might be more likely to be satisfied than the ones above.
For example, if someone introduces you as a mathematician, then you will be given the identity of mathematician.
In this case, this identity might be ephemeral, lasting a mere conversation.

\subsection{Generating Further Definitions}
In the above sections, we have presented three different axes under which to establish who is (and who is not) a mathematician. 
Much like the generating set of a group that can build other elements, each of the definitions we have given above could lead to further definitions about whether someone is a mathematician or not.
Most individuals use a combination of these definitions to self-define or define others as ``mathematicians." 
In this subsection, we will now consider the previous three definitions as a generating set, then look to see if some union, intersection, or other combination is satisfied.
Because this new definition refers to sets generated from the other mathematician definitions, we refer to this new generated definitions as \emph{compuesto} definitions. 

\newtheorem{compuesto_ex_definition}{\emph{Compuesto Definition}}
\begin{compuesto_ex_definition}[Union/broad]\label{def:compuesto-broad}
    A mathematician is someone who satisfies any definition of a mathematician. Effectively, this definition is the union of all definitions.\\
    
\end{compuesto_ex_definition}

\begin{compuesto_ex_definition}[Intersection/narrow]\label{def:compuesto-narrow}
    A mathematician is someone who satisfies every definition of a mathematician. Effectively, this definition is the intersection of all definitions.\\
\end{compuesto_ex_definition}
\noindent It is entirely possible that this narrow definition would result in the empty set. 
We investigate this possibility by examining the identities of the co-authors through their positionality statements in Section~\ref{subsec:positionality}, with the results summarized in Table~\ref{table:position}.

\section{Analysis}
\label{sec:analysis}

In this section, we analyze the example definitions provided above to explore the similarities and to tease apart the differences between both the {\em types} of example definition (function-based, qualification-based, and identity-based) as well as the dichotomous variation available for each one (i.e., either narrow/broad or internal/external). 
Our analysis takes two forms. Our analysis moves from the micro to the macroscale, first by applying these definitions to individual co-authors of this paper and then moving to the larger discipline-specific implications of these example definitions. 
This section is organized as follows. First, in Section~\ref{subsec:positionality}, we offer positionality statements for each of this paper's co-authors and then situate their relationships to each example definition type (and axis). 
Then, in Section~\ref{subsec:implications}, we use a polarity mapping analysis to explore implications of using the example definitions from Section~\ref{sec:defs} when defining a ``mathematician''. 
Specifically, we used a participatory process to capture the authors' perspectives on the opportunities and risks associated with the application of extremal versions of each definition type.
Finally in Section~\ref{sec:AndBeyond}, we bridge from this broadening analysis of the example definitions (Section~\ref{sec:defs}) towards future conversations and implications on the whole field of mathematics. 

\subsection{Positionality Statement}\label{subsec:positionality}

Here, we offer a series of personal examples of self-definition via individual positionality statements to help readers contextualize our work. 
We then use the definitions from Section~\ref{sec:defs} to analyze these positionality statements and provide a table that summarizes this analysis. 

\subsubsection*{Ron Buckmire}
I am a cisgender, gay man of Afro-Caribbean descent who has been a member of the mathematics faculty of a small liberal arts college in Los Angeles, California since 1994. 
I earned all three of my degrees in mathematics at Rensselaer Polytechnic Institute, an Engineering-oriented institution.
Perhaps influenced by where I was educated, when I do mathematics, I often think about it in the context of solving a specific problem, and almost always use computation to assist me in the solution process.
I believe that “mathematics is for solving problems” and that ``mathematics is a human endeavor." 
Mathematics is done, taught, learned, discovered and created by people, and that the identities and prior experiences of people who do math, are important.

\subsubsection*{Carrie Diaz Eaton}
I am a mathematician by degree. 
Each higher degree I hold is in mathematics, but all my degrees were very interdisciplinary with biology.
I am not a mathematician by job title, having had appointments in biology and computational studies departments after graduation. 
I used to joke that I was the biologist in the room of mathematicians and the mathematician in the room of biologists, so my disciplinary identity feels fluid and relative at times, much like my queer and multi-ethnic identities, respectively.
My primary research activity is in STEM education. 
My methods are statistical, computational, and qualitative. 
I do not partake in formal proofs in any of my scholarship, yet all of my work requires sound logic.
I am deeply involved in my mathematics community, but also involved in the biology education community as well. 
As a Latina, I resonate with the idea of Nepantla \cite{gutierrez2017} and see boundary spanning as my superpower. If I relied on all definitions of a mathematician or any discipline to be satisfied, then I would have no disciplinary identity.
Instead I choose to embrace all my boundary-spanning disciplinary identities as a strength.

\subsubsection*{Joseph E. Hibdon, Jr.} 
I am a cisgender Native American man, who is the descendant of the Luiseño Band of Mission Indians. 
My Ph.D. is in Applied Mathematics.  
On a daily basis I teach and do research in applied mathematics at a public, four-year, university that is a Hispanic Serving Institution. 
Much of my research is in the realm of mathematical education and I do not do proofs; beyond the courses I am teaching. 
In general, I identify myself to others as a practicing mathematician who works to make mathematics available to all.

\subsubsection*{Katherine M. Kinnaird}
On a daily basis, I think about computational tasks.
This can include: scoping a project to leverage existing (or soon-to-be-collected) data, drafting an algorithm, creating a method for collecting data, coding an algorithm, analyzing data, or writing up results.
My bachelors, masters, and Ph.D.~are in mathematics. 
In my Ph.D., I did prove a lemma and a few corollaries, but I have yet to publish a proof in a peer-reviewed journal. 
While others might call me a mathematician due to my qualifications, I would no longer call myself a mathematician.
I identify as a computational scientist, by which I mean a person who exists in the intersection of data science and computer science. 


\subsubsection*{Drew Lewis}
I am a cishet white male mathematician. 
My early career work was squarely in the realm of ``pure'' mathematics (algebraic geometry), but I am now probably better described as a mathematician who does education research. 
This latter work has shown me how who is perceived as a ``mathematician'' has real implications, as often the views of mathematicians (including me) are given disproportionate weight in math education spaces.

\subsubsection*{Jessica Libertini} 
My day-to-day activities include overseeing interdisciplinary research efforts that use creativity science, design thinking, and operational planning to address global challenges. 
While these activities rarely are ``math'' in the strictest sense, I often draw on mathematical, modeling, computational, and logic frameworks to advance progress.
I do not write proofs and have limited experience with proofs, as I am trained in mechanical engineering (B.S. and M.S.), applied mathematics (Sc.M. and Ph.D.), and international relations (MAS).
I do not consider myself a mathematician; rather, I identify as an undisciplinarian who challenges disciplinary boundaries.

\subsubsection*{Omayra Ortega}
My day-to-day activities as an applied mathematician include teaching mathematics and statistics courses, conducting research with my collaborators at all levels (undergraduates, grad students, post-docs, faculty, and other professionals). 
I rarely write proofs unless it is in the classroom setting. My research focuses on solving real world problems and relies heavily on data analysis and model fitting. I have a Ph.D. in Applied Mathematics \& Computational Sciences and degrees in music, ``pure" math, biostatistics and public health. Broadly, I identify as an Afro-Latinx non-binary women of Panamanian descent with strong New York roots.

\subsubsection*{Rachel Roca}

On a day-to-day basis, outside of my course work, I spend much of my time reading and discussing articles (mostly from applied math and STEM education), working with code, and TAing a computational methods course.
My work, thus far, has been focused primarily on course work, some of which are proof-based, while my current research is more computational based. 
I graduated with a bachelors in mathematics from a small liberal arts school in 2021 while also minoring in computer science and Spanish.
Currently, I am a second year Ph.D. student in computational mathematics, science, and engineering at Michigan State University. My research interests span topological data analysis, computing education, and math for social justice.
In a broader sense, I identify as a queer white woman mathematician (in training).
Within my department, I am a mathematician, but do not feel like one among those in a traditional math department. 
Due to being a newcomer to the world of academia, coupled with my interdisciplinary work, I am still in the process of exploring and developing my own identity.

\subsubsection*{Andr\'es R. Vindas-Mel\'endez}
I am a queer, chronically-ill, Costa Rican-American mathematician raised in South East Los Ángeles, California.
I am a first-generation college graduate having, earned my bachelor's degree in mathematics and minored in Philosophy and Chicanx/Latinx Studies. 
I also earned a master's and Ph.D. in mathematics.
While pursuing a Ph.D, I earned a graduate certificate in Latin American, Caribbean, \& Latinx Studies. 
All of my schooling has been through public institutions. 
My research is primarily in combinatorics and I have been expanding my research interests to include applications of data science to combinatorics and mathematics for social justice. 
In particular, for my combinatorial work I prove results on lattice-point enumeration and triangulations of lattice and rational polytopes.
In mathematics for social justice, I am interested in data analysis, interdisciplinary study (e.g., social science, history, economics), and development of mathematical techniques to investigate racial/social issues.

I believe mathematics is as diverse and dynamic as those who use it.
My mathematics is influenced by my experiences and I bring to it my own perceptions.
I aim to build meaningful and empowering experiences with mathematics, while also challenging others to think about the power structures that are present in and outside mathematical spaces. 

\subsubsection{Reflections on the Co-Authors' Positionality Statements}

We examined our positionality statements through the lenses of each of the definitions outlined in Section \ref{sec:defs}. 
When we analyzed the positionality statements given above, written by each of this article's authors, and collapsed the findings into a table, we found some interesting patterns. 
We would expect this group of co-authors to have many similarities, since we began this work as research fellows at the Institute for Computational and Experimental Research in Mathematics (ICERM), a prestigious mathematics institute, in Summer 2022.
In fact, all of the co-authors were classified as a ``mathematician" under the following seven example definitions: Functional-Broad, Functional-Hybrid, Qualifications-Broad, Qualifications-Hybrid, Identity-External, Identity-Hybrid, and \emph{Compuesto}-Broad; therefore those columns are not included in Table \ref{table:position}.
In these seven common definitions there are both kernels of what it generally means to be a mathematician and qualities unique to individuals who have been entrenched in, and validated by this field for many years.
We include Table \ref{table:position} to highlight the diversity within our group while, by omission, also highlighting the commonalities. 
\begin{table}[ht]
\begin{tabular}{p{0.8in}p{0.8in}p{0.8in}p{0.8in}p{0.7in}}
\hline\hline
Author & Functional-Narrow & Qualifications-Narrow & Identity-Internal & \emph{Compuesto}-Narrow\\
\hline
RB & & X & X & \\
CDE & & X & X & \\
JEHJ & & X & X & \\\hline
KMK & & X & & \\
DL & X & X & X & X\\
JL & & & & \\\hline
OO & & X & X & \\
RR & X & & X & \\
ARVM & X & X & X & X \\
\hline
\end{tabular}
\label{table:position}
\caption{This table highlights the differences among the authors alignment with each Example Definition.}
\end{table}

Many of us identified at one time as applied mathematicians, who are often excluded by narrow function-based definitions due to a lack of proof writing, and we explicitly mentioned this in our positionality statements. While, intellectually, we know that one narrow functional definition does not need to define us, as a group, many of us leaned into other definitions of mathematician. 
Perhaps we were just speaking to the myriad of definitions described, but for some of us, there is trauma around whether external entities would identify us as mathematicians, particularly when we hold other marginalized social identities in mathematics.

When narrow functional definitions were not satisfied, many suggested using broader functional definitions (e.g., RB and DL), either explicitly or implicitly through describing how they use mathematics (e.g., JEH, CDE, and OO). 
Most co-authors also leveraged other definitions of mathematician such as qualification-based and identity-based (see Table \ref{table:position}). 
RB, JH, and AVRM
describe beliefs and future hopes that mathematics should be accessible for all or should used in service of humanity.
Redefining boundaries through futurism is a revered tradition among scholars who study marginalized identities and power (\textit{i.e.}\cite{anzaldua1987borderlands, maldonado2014futuristic}).

Those of us with pure mathematics training addressed power differences based on identity in other ways by thinking about power between mathematics and mathematics education or in the field of mathematics more generally. 
For example, CDE used to say ``I am a mathematician in a room of biologists and a biologist in a room of mathematicians." 
This is an example of identity as relative to the context, if alternatives are available. 

This same quote also invokes another observation on identity as expressed in the positionality statements -- that they are fluid over time. 
Early career co-authors readily acknowledged that their explorations of their identities as mathematicians are ongoing. 
KMK, a mid-career co-author, expressed leaving a mathematician identity for that of a computational scientist, reflecting an evolution in research area over her career.

Finally, we see some co-authors willingly accept multiple simultaneous disciplinary identities. 
This could be expressed as a subdisciplinary identity, such as ARVM describes, or could be from different disciplines, such as CDE describes.
Alternatively, JL feels limited and constrained by disciplinary labels and therefore rejects any disciplinary identity in favor of ``undisciplinarian.''

Another challenge to the analysis of the positionality statements is the difference between the definition of a mathematician that an individual may hold personally, versus the definitions that an individual may feel beholden to from outside. 
Many acknowledged that their positionality statement listed the qualifications that they believed most in the mathematical sciences would accept as the correct qualifications or qualities for someone to be considered a mathematician.  
This internal versus external duality adds a level of complexity to the analysis of an individual's positionality statement.

\subsection{Implications of Definitions: A Polarity Analysis}\label{subsec:implications}

During the authors' time at ICERM in Summer 2022, we used a polarity framework to help us explore the tension between inclusive (broad or internal) and exclusive (narrow or external) definitions of the term ``mathematician."
This framework was introduced by Barry Johnson in the 1990s and has been widely adopted in the business leadership and design thinking communities \cite{polarities}. 
A polarity analysis is a participatory methodology that allows a group to explore the tensions that exist when there is value and wisdom in two sides of seemingly polar opposite perspectives or approaches. 
In this case, we used the polarity analysis to understand how those of us working to advance social justice in mathematics might find value in all of the definitions, depending on our purpose.

The author team's approach was very similar to the Polarity Mapping exercise outlined by University Innovation Fellows (\url{https://universityinnovation.org/wiki/Resource:Polarity_Mapping}).
First we defined each pole; in this case the two poles are the broad definition(s) of mathematician on one end and the narrow definition(s) of mathematician on the other. 
(Note that in the case of identity-based definitions, we classify internal as broad, as the autonomy is given to the individual, and external as  narrow, as the power to bestow the title relies on others who may apply a definition that would exclude a self-identifying mathematician.)

Then we spent time individually ideating on the opportunities of focusing on each pole as well as the risks of over-focusing on a pole at the expense of energizing the other pole. 
Next we discussed the results of the individual ideation in small teams, and ultimately, we shared these with a larger group that included our co-author team as well as other participants in the ICERM program. 
Each idea was captured on a sticky note (see figure \ref{fig:polarity} in the appendix), and this allowed the team both to appreciate areas of high synergy (where multiple team members had identified a common theme) as well as to highlight those creative ideas at the margins that may have only been thought of by one individual. 
We then discussed what had been revealed in this initial analysis as an entire group.

Lastly, the team also explored action steps that could be used to energize the positives of one pole, as well as early warning signs that a pole could be over-emphasized at the expense of the other pole. 
This portion of our polarity analysis was not as well-suited to our efforts, as we are not beholden to adopting a single definition; in fact, we argue for intentionality in selecting, stating, and justifying a definition unique to the purpose at hand. Therefore, the action steps and early warning analysis has been excluded from this paper.

Based on our polarity exercise as well as our work in exploring different types of definitions, we specify potential motivations and the implications of using narrow or broad definitions through the lens of opportunities and risks below. We acknowledge that some team members entered the exercise with skepticism, voicing a strong preference for more inclusive definitions.  However, the analysis revealed opportunities and risks for each type of definition.
We also remark that due to the participatory nature of this qualitative research approach, the results reflect the perceptions, beliefs, and values of the participants.
Since this work was done in the context of a research program on social justice and mathematical sciences, the co-authors' brought a social justice lens to the polarity mapping exercise and its analysis.

\subsubsection{Opportunities of Narrow Definitions}
\label{subsec:positive-narrow}

A narrow definition may initially seem limiting, but there are purposes and settings that justify using a more narrow definition to support social justice efforts. 
Below, we explore each narrow definition and its utility towards social justice.

The narrow Function-based Definition~\ref{def:functional-narrow}  (``A mathematician is a person who proves theorems using proof techniques") allows for there to be a well-defined cutoff for membership in the class of mathematicians, so that measurements, data collection, and statistics are conducted in a relatively narrow context. 
Such studies could leverage the narrowness of the definition to highlight gender, racial, and other representational disparities, study the impacts of interventions designed to address the disparities, and compare the disparities in the mathematical community with those in other narrowly defined professions such as medical doctor or trial lawyer.

The narrow qualifications-based Definition~\ref{def:qualification-narrow} (``A mathematician is someone who has a
Ph.D. in mathematics and currently holds a research position in mathematics.") also has the benefit of having a well-defined boundary that could be used for social science research about the field of mathematics.
For example, one could use this definition to conduct surveys or interviews with individuals that have a Ph.D. in mathematics and also do active research in a defined mathematical field.
Such a study could help reveal insights about how welcoming (or unwelcoming) the narrowly-defined community is to persons from underrepresented communities who meet the membership criteria; a similar study could explore the differences in experiences between Ph.D. research mathematicians who are male and white with those who are not.

Similarly, from a justice-oriented perspective, using narrow function-based and qualification-based definitions can allow us to more easily identify who is being excluded from the mathematics community.
Using a qualification-based definition of earning a Ph.D. in mathematics, it is widely documented that woman and people of color are underrepresented in this version of ``the mathematics community" as compared to white men \cite{Buckmire2021}. 
As the mathematics community is heterogeneous with many disciplinary areas contained within it, we might also consider what kinds of individuals are most likely to be excluded from different subareas of mathematics when using a narrow function-based definition of ``mathematician" (like ``a mathematician is someone who publishes papers with proofs of theorems").
 
 For the identity-based definitions, one could consider the external definition to be the equivalent analogue of the narrow modality for the function-based and qualifications-based definitions.
 Recall, the external identity-based definition, \ref{def:identity-external}, is ``A mathematician is anyone whom others say
is a mathematician."
 We found it difficult to develop a social justice benefit to having groups excluded in a potentially arbitrary way from the mathematical community.
 However, we noted that there was at least the potential for transparency with an opportunity to gain insights into the perspectives of those who draw a tighter circle around who they consider to be a mathematician.
 For example, in a paper, the definition used by the authors may be clearly stated.
 Alternatively, even in casual conversations, insights and contexts can be gleaned to reveal the dynamics of social perception about who counts as a mathematician by different individuals and groups.
 As we participate in formal and informal conversations about who is a mathematician, those working to promote social justice within the field of mathematics can monitor the social dynamics in various ecosystems throughout the mathematical community.
 This has the potential to illuminate which barriers are eroding and which are not, allowing more focused efforts in distributing the power of mathematics.
 
This subsection has demonstrated that there can be social justice benefits to adopting a narrow definition of the term ``mathematician" under certain conditions, despite the exclusive tone that is inherent in a narrow definition.
Specifically, narrow function-based and qualification-based definitions offer clarity for researchers seeking to shine a light on inequities.
With the potential for systemic problems to be identified through the use narrow definitions, they may lead to the creation of a more equitable, just mathematics community.

\subsubsection{Risks of Narrow Definitions}
\label{subsec:negative-narrow}

As shown above, there are some benefits to the adoption of a narrow definition of the term ``mathematician," but there are also risks.
Below we explore social justice risks of employing the extremal (narrow and external) function-based, qualification-based, and identity-based definitions of a mathematician.

The narrow Function-based Definition~\ref{def:functional-narrow} (``A mathematician is a person who proves theorems using proof techniques.") does not allow for an inclusive interpretation of the mathematics community beyond some subjective definition of what a mathematician ``should' be able to do (i.e., prove theorems).
Those that may not have had the privilege of having a proof-based course or are just getting started in the exploration of mathematics would not be considered \textit{part of the mathematics community} under Function-based Definition~\ref{def:functional-narrow}.
It is likely that the people who are excluded from the community through the application of Function-based Definition~\ref{def:functional-narrow} are more likely to be people of color, women, and people who have less access to proof-based courses and may have attended lower-resourced institutions. 
Thus, the narrow function-based definition would generally tend to reduce the diversity (and maintain the overrepresentation of maleness and whiteness) in a community defined by this definition. 
This particular effect is a substantial negative outcome of the use of a narrow function-based definition.

Another risk of the use of the narrow Function-based Definition \ref{def:functional-narrow} is that this definition fails to include the students and graduates that increasingly populate data science and other programs that are very applied, but may not have a proofs requirement.
By exclusion of these individuals, we also lose the opportunity to count (as mathematics) the work that they do, which has the potential to contribute more tangibly to the social fabric, as data scientists are explicitly trained to use data to understand the world around us.

Examining the narrow Qualification-based Definition \ref{def:qualification-narrow} (``A mathematician is someone who has a Ph.D.in mathematics and is currently doing research in mathematics") also exposes several risks to the advancement of social justice in mathematics because it imposes barriers on who is part of the in-group and leaves many in the out-group.
If, for example, we employ an alternative narrow qualification-based definition, ``A mathematician is someone who has earned a Ph.D. in Mathematics at an Ivy League institution," then those who meet the criteria may feel a strong community bond, which can exacerbate the sense of exclusion for those who do not belong.
As another example, research in faculty hiring in mathematics (and other disciplines) demonstrates a cycle of elitism and exclusion through the application of a narrow qualification-based definition, such as ``A mathematician (worthy of being hired at a prestigious university) is one who currently holds a mathematical research position at a(nother) prestigious university."
The result is that high prestige universities exchange faculty with one another or export faculty to less prestigious schools but rarely hire faculty from the less elite universities \cite{Wapman2022}.
Furthermore, this definition may dissuade future students from selecting a mathematics major, as this definition of mathematician sharply limits what careers they could pursue, despite the data that shows mathematics majors being successful in a broad range of career choices beyond those meeting the narrow qualification-based definition \ref{def:qualification-narrow}.

In analyzing the external Identity-based Definition \ref{def:identity-external} (``A mathematician is anyone whom others say is a mathematician"), we reveal risks similar to those of the \ref{def:qualification-narrow}.
These definitions require an external review, and unlike the function-based definition, this review of who qualifies to be a mathematician can lead to subjective barriers of belonging in the mathematical community, which can breed elitist and exclusive mentalities that marginalize those who do not meet the criteria.
This means that these types of narrow/external definitions, if employed without care, could perpetuate the marginalization of those who do not satisfy the external criteria of a selected definition.

The negative implications of a narrow definition are not merely static, as they provide a feedback loop that can exacerbate exclusivity. 
For example, leveraging the theory of possible selves, potential future members of the community may self-select out if they do not see people with identities similar to their own in the in-group \cite{Oyserman}.

\subsubsection{Opportunities of Broad Definitions}
\label{subsec:positive-broad}

Here, we will discuss the social justice opportunities that emerged through our polarity analysis as we examined function-based, qualification-based, and identity-based broad definitions definitions of ``who is a mathematician."

Recall our broad Function-based Definition~\ref{def:functional-broad}, ``A mathematician is a person who uses mathematical concepts, tools or techniques to study and solve problems."
Such a definition can foster a mathematical community that is welcoming by including and engaging those who are otherwise excluded from the conversation of ``who is a mathematician." 
When we apply a broad function-based definition, we see that the mathematical community that results now includes those who are not necessarily writing proofs but are using mathematical tools in contextual ways -- as data scientists, statisticians, and quantitative social scientists.
The broader the circle, the more people are included, which provides an opportunity to put the power of mathematics (and the power of the title ``mathematician") into the hands of people with a collectively broader perspective than the one dominated by a smaller (and largely white and male) community.
This approach also allows for broader recruitment and engagement and could support a migration from an environment where certain folks are relegated to disciplinary silos to a more collaborative interdisciplinary environment.

By exploring Broad Qualification-based Definition \ref{def:qualification-broad} (``A mathematician is someone who has completed a course in mathematics at any level"), we see that this definition is highly inclusive.
This level inclusion can be used to empower people of all ages, but it has the potential to be particularly useful in planting seeds of possibilities into the minds of young children who are just starting to consider their future careers.
Using a broad qualification-based definition offers the opportunity to rebrand mathematics as something that everyone can do, invest in a growth mindset, and challenge the societal norm that ``math is hard."
By investing in the mindsets and perceptions of young children and others who have not previously considered themselves to be mathematicians and inviting them into the circle of mathematicians, there is the opportunity to begin dismantling some of the barriers that exclude people from seeing their future self as a mathematician - because they have been told they already are one.
Such an approach has the potential to challenge the hegemony of white males in mathematics community that can be created through the application of a broad qualification-based definition.

When considering a ``broad" identity-based definition, we selected the Internal Identity-based Definition \ref{def:identity-internal}, ``A mathematician is anyone who says they are a mathematician."
This definition inherently empowers individuals to opt-into the community of mathematicians, and this has the potential to increase the diversity of that community.
Recall that one of the risks of narrow definitions in general was that of driving away potential community members who do not see themselves represented in the group. Conversely, a broad definition casts a wider net and therefore paints a picture of ``who is a mathematician" that can be more diverse. 
And that increased diversity has the potential to be more inviting to more members who may have otherwise opted out of engaging with the community when it was defined narrowly. 
The people who would be most affected by such a change include persons of color, women, and people from other underrepresented identities. 
Note that this concept is the contrapositive to the point on possible selves and the non-static impacts of an exclusive or narrow definition given in the previous section \ref{subsec:negative-narrow}. 
So, in this case broadening the definition of mathematician makes it more likely that the community will continue to grow in diversity.

\subsubsection{Risks of Broad Definitions}
There are also downsides to the application of the broad version of the definitions of ``who is a mathematician."
In this section, we explore the downsides to broad function-based, qualification-based, and identity-based definitions of ``who is a mathematician."

As previously noted, the Broad Function-based Definition~\ref{def:functional-broad} (``A mathematician is a person who uses mathematical concepts, tools or techniques to study and solve problems") is very inclusive by design.
However, it could be so inclusive as to include nearly everyone to the point of dilution of the value of the term at all.
For example, is a farmer who is estimating their yield a mathematician? 
Using a broader definition of mathematician may dilute the power of those who enjoy this access under the current definition; this  potential loss of power may explain why certain groups may cling to narrower definitions of the term and oppose broader ones. 
On the other hand, using a broader definition could a rhetorical choice to intentionally dismantle this power.

For this reason, it may be appropriate to discuss the Hybrid Function-based Definition \ref{def:functional-hybrid}, ``A mathematician is a person who, as part of their daily work, employs mathematical techniques and tools to solve mathematical and other problems."
However, even this definition could be so inclusive that it may result in people who do very different kinds of mathematics (i.e., who use very theoretical and very applied techniques) being lumped together under the same umbrella.
When people who view themselves as being very different are given the same label face practical challenges when their internal sense of who they are is incongruent with the external impression of how they are viewed by others, such as more limited access to resources (\cite{jing09}) and more time taken to cultivate external identity (\cite{lam2020hybrids}). 
The negative consequences of an identity crisis fueled by use of a broad function-based definition could include: lack of confidence to do mathematics, reduced identification with the mathematics community and increased fragmentation into disciplinary silos.

The broad Qualification-based Definition~\ref{def:qualification-broad} (``A mathematician is someone who has completed a course in mathematics at any level") also has risks.
If the definition of mathematician is broad to the point of near universality through the application of the Internal Identity-based Definition~\ref{def:identity-internal} (``A mathematician is anyone who says they are a mathematician"), then the barriers that have excluded people, particularly those who are not white and male, are artificially erased. 
The barriers to entry and the pressures to leave the mathematics community may still be very real and marginalizing to some, so the adoption of such a broad definition risks erasing the perception of those barriers and pressures, reducing access to resources to dismantle those barriers and pressures, and minimizing the acknowledgment of harm to those affected by those barriers and pressures. 
In other words, while such a broad definition might work in a utopia, there is risk of exacerbating harm if applied in our current world.  

The Internal (Broad) Identity-based Definition \ref{def:identity-internal} allows anyone to self-identify as a mathematician.
While this definition supports individual agency, it also risks the dilution of a sense of community in mathematics as it is difficult to assess the commonalities of the resulting body of people.
This also could impact the resources available to the community, as there is a risk of less homogeneity in the definition of who a mathematician is, what a mathematician does, and therefore, what a mathematician needs to be successful in their pursuits. 
What is more, this kind of dilution is reminiscent of the risks that Guti\'errez\cite{gutierrez} discusses in terms of dilution of power. 

\subsubsection{Overview of Polarity Analysis}

\begin{figure}[ht]
    \centering
    \includegraphics[width=1.0\textwidth]{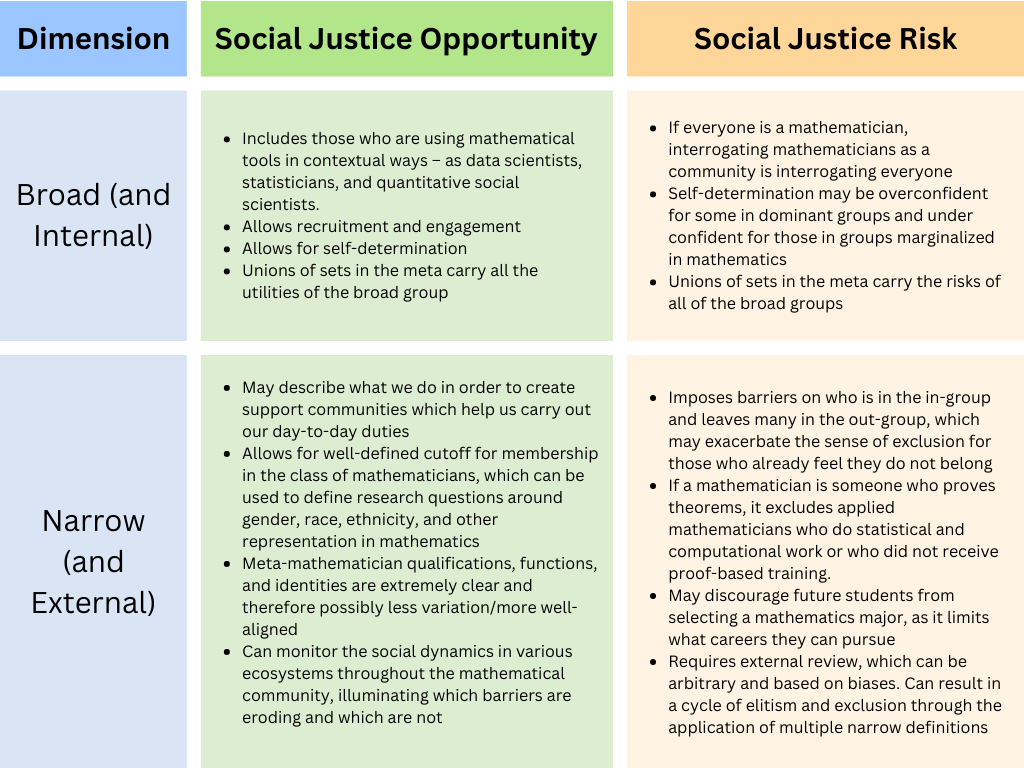}
    \caption{A summary of some of the opportunities and risks associated with broad, narrow, internal, and external dimensions of the definitions proposed.}
    \label{fig:summary}
\end{figure}

The polarity analysis revealed wisdom and opportunities, but also risks as we explored the extrema of the continuum of definitions along the narrow-broad (external-internal) continuum. 
The overall findings are collectively captured in Figure \ref{fig:summary}.
As we look across these findings, we note some key themes.

In the pursuit of social justice, narrow definitions provide those working towards social justice in the field an opportunity to perform research that illuminates systemic barriers and pressures that disproportionately exclude certain demographics from inclusion under those narrow definitions.
However, they also carry the risk of implying that ``real mathematicians are ...", causing new harm and exacerbating existing harm.
Therefore, narrow definitions should be used with caution and with clear explanation about why they are being employed.

Conversely, broad definitions, by their very nature, are more inclusive, which invites more diversity into the community and provides a positive feedback loop welcoming in future members who see themselves represented.
However, these definitions could lead to a dilution of identity cohesion in the mathematical community (loss of an understanding of what/who is a mathematician), reduce available resources, and risk fragmenting the community into smaller, more homogeneous groups.
Perhaps even more harmfully, if used without care, broad definitions could imply that the systemic barriers and pressures that have historically been used as gatekeepers against certain groups (e.g., women and persons of color) have simply disappeared because, well, now anyone can be a mathematician.
This could actually make things worse, as it could minimize the lived experiences of those encountering marginalization.
Therefore, broad definitions should also be used with caution and with clear acknowledgement that some may continue to face barriers.

We also see that there is a benefit to using hybrid definitions. 
For example, when we explored the functional definition, the narrow definition carried the risk of excluding those working in mathematics beyond the world of proofs, but the broad definition was so wide-ranging it risked defining someone who uses mathematics to calculate the tip on a restaurant check as a mathematician.
If there is a need for a cutoff for the definition of ``a mathematician," but one needs to, simultaneously, cast a large net, then the Hybrid Function-based Example Definition~\ref{def:functional-hybrid} (``A mathematician is a person who as part of their daily work, employs mathematical techniques and tools to solve mathematical and other problems.") is one appropriate definition that could be employed. 
Many people getting degrees in mathematics and related fields will not go on to become professors of mathematics and will not spend a career writing proofs, but they \textit{will} be using mathematical skills in their careers.  
Having students realize that, even if they pursue careers outside of academia, they will still be part of the mathematical community has significant benefits in that it enlists many more people to the profession and thus solidifies the future of the mathematical community.
At the same time, drawing a boundary has the potential to create a sense of community for this larger group of those who meet this hybrid function-based definition.

Lastly, we note that the application of any definition of a mathematician is shaped by the motivation behind its selection. 
For example, in this article we are engaging in an open conversation about defining the term  ``mathematician" as a means to encourage others to join the mathematical community.  
But people with other motives could use these definitions to accomplish different ends. 
We therefore advise those involved in this discourse to be mindful of the implications of their own definition choices and to use a critical lens when examining the definitions of others.

\subsection{From Definitions to Future Impacts}\label{sec:AndBeyond}

Above we have given examples of how various definitions may be used both by individuals and applied to larger groups of people. 
As with all mathematical definitions, these definitions reflect choices made by people, and we have explored some of the implications of those choices.
Importantly, a choice of definition depends heavily on context and the motivations of the person/people choosing that definition. 
Thus, when presented with a definition, it is important to interrogate: who is making the definition?  In what ways does this choice work towards or against justice?  What groups benefit from the choice being made, and who is harmed?

In addition to asking readers to apply a critical lens to how mathematician is defined, by themselves and by others, the authors encourage readers to engage in self-reflection as well as discussions within their own communities - with students, colleagues, or at mathematics conferences.
Perhaps set aside time for yourself (and others) to explore the following series of questions:
\begin{itemize}
    \item What is my positionality statement? Do I consider myself a mathematician, and how is the answer to this question important to my own identity?
    \item Using a function-based definition, do I meet the narrow definition? If not, how might I develop a hybrid definition that captures what I do? And under that new definition, who is included, and who is excluded?
    \item Using a qualification-based definition, do I meet the narrow definition? If not, how might I develop a hybrid definition that captures what I do? And under that new definition, who is included, and who is excluded?
    \item Using an identity-based definition, who would define me as a mathematician, in what settings or situations, and why?
    \item How do I define mathematician? How does my institution, employer, and/or social group define mathematician? Who is included and who is excluded with each of these definitions, and what does that mean both for the excluded individuals and for the resulting community?
    \item As we see a surge in programs in data science and other similar fields, how might intentional decisions about including or excluding these groups in the definition of mathematician shape the mathematical community?
\end{itemize}

\section{Conclusion: Identity as Power}
\label{sec:conclusions}
Mathematics and the epistemology of mathematics is socio-cultural by nature, rather than absolutist, argues Burton \cite{BurtonFem}. 
As we prepared to move into research projects in which we would collect and analyze data on mathematicians, we were forced to try to make assumptions about who is a mathematician, in order to provide guidance in determining whose data, which data and what data was worth collecting for our various social justice and data science projects. 
Our search for the right definition of a mathematician eventually came to an unsatisfying, but possibly predictable answer: it depends. 
Your choice of definition for mathematician will depend on why you are collecting the data and which questions you are answering. 
It will depend on what data is available to you for collection purposes.  
It will depend on whether your purpose is to use the definition to exclude individuals from opportunities, to welcome people into a more inclusive spaces for all mathematical identities, or to illuminate the barriers that have make it difficult for some individuals and groups to become and persist as mathematicians.

In determining who counts as a mathematician, Burton highlights the power of Eurocentric male dominance in the mathematics community, as well as the exploitation and erasure of those who do not share those identities. 
Those who are ``mathematicians'' wield the power to define and identify ``important'' mathematical areas, where value is accorded to some results rather than others and decisions are taken on what should or should not be published in a society determined by ''power relationships'' \cite{BurtonFem}.
In the context of our own research group, in order to create the field of data science for social justice within mathematics, we had to assert ourselves as mathematicians (and a computational scientist and an undisciplinarian) to define the field. 
Therefore, the rationale for determining who is a mathematician moved beyond one for data research purposes and became directly connected to defining the worth of our research program.

All of these decisions and mathematical contributions are products of people, who may include or exclude others from the mathematics community. 
Possible exclusion from the traditionally elite math community, what Bartholomew et al.\cite{Bartetal} call the ``maths club," threaten identity formation. 
This ``maths club," utilizes fear and narrow definitions of mathematicians to keep power to some and disempower others.
Even advanced undergraduate mathematics students express uncertainty about their identity in being a mathematician \cite{Bartetal}.
In addition, research has repeatedly demonstrated that co-construction of social identity and mathematics identities are tightly related. For example, a case study of a student's experiences navigating and overcoming racialized treatment that aimed to exclude her from the mathematics community by Oppland-Cordell \cite{OpplandCordell}, highlights how racial, gender, and class identities affect mathematics participation and learning. 
Therefore, we want to bring attention to how crafting mathematical environments that embrace and showcase the strengths of multiple social identities is just as important as creating environments that nurture a variety of mathematical identities. 

\section{Future Work}
\label{sec:future}

In the preceding, we have explored various definitions of ``mathematician'' and analyzed the implications of these kinds of definitions. 
While our definitions emerged from the perspectives of a relatively diverse group of mathematicians, we are only a small subset of the mathematical community. 
As such, we anticipate this work as a first step for future work exploring how the broader mathematics community interprets the term.  
Some questions we anticipate answering are: 
\begin{enumerate*}[label=\arabic*)]
\item To what extent does the mathematics community have a shared definition of ``mathematician''?  
\item Do members of the mathematics community switch among multiple definitions to suit particular contexts, or persist with a rigid idea of who is a mathematician? 
\end{enumerate*}
As discussed above, different answers to these two questions have different consequences for the profession and the people who are or would be part of it.

What is more, we hope this work will inspire the mathematical community to think both locally and globally about the definitions of ``mathematician'' we assume and use, and the consequences these definitions can have in the broader community as well as in our classrooms, majors, and institutions. 
With this in mind, we invite readers join us in examining identity in mathematics either individually (in one's sphere of influence) or collectively (perhaps in collaboration with the authors of this paper). 
Mathematicians, as those in other disciplines, often accept the terms of their profession without question or thought. 
We hope this paper has inspired you to question these assumptions. 
If we are more thoughtful about how we define ourselves, our surroundings, and the community that we work in, we can influence the social contract of this community. 
Put another way, this is a call to action to create the mathematics community that you wish to see.
If you are interested in joining our working group, please reach out to the authors on this paper.

\section*{Authorship Contributions}

All the listed co-authors contributed equally to all portions of the paper.

\section*{Acknowledgements}
This work was supported in part by the National Science Foundation under Grant No. DMS-1929284 while the authors were in residence at the Institute for Computational and Experimental Research in Mathematics in Providence, RI, during the Data Science and Social Justice: Networks, Policy, and Education  program.
ARVM is partially supported by the National Science Foundation under Award DMS-2102921.
RB acknowledges the support of the Occidental College Office of the Dean of the College for funding his sabbatical in the 2022-2023 academic year.
KMK is the Clare Boothe Luce Assistant Professor of Computer Science and Statistical \& Data Sciences at Smith College, this work has also been partially supported by Henry Luce Foundation's Clare Boothe Luce Program.

\bibliographystyle{plain}
\bibliography{bibliography}

\section*{Appendix}
\begin{figure}[h]
\begin{center}
    \includegraphics[scale=0.08]{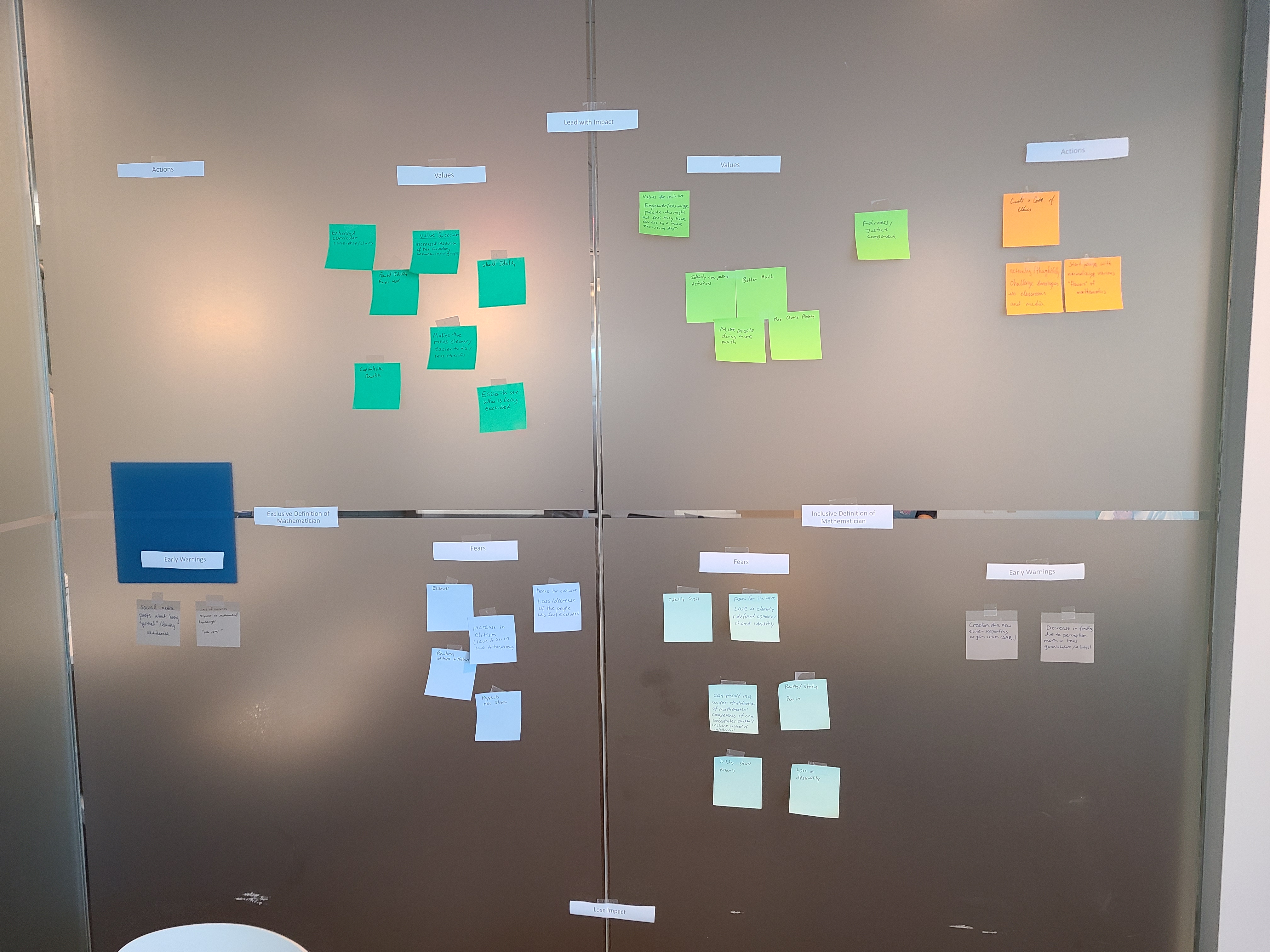}
\caption{Documentation of the polarity analysis conducted by the authors.} \label{fig:polarity}
\end{center}
\end{figure}
\end{document}